\documentclass[11pt]{amsart}
\usepackage{amsmath}
\usepackage{amssymb}
\usepackage{amscd}
\newtheorem{dfn}{Definition}[section]
\newtheorem{rem}[dfn]{Remark}

\newtheorem{thm}[dfn]{Theorem}
\newtheorem{lem}[dfn]{Lemma}
\newtheorem{prop}[dfn]{Proposition}
\newtheorem{cor}[dfn]{Corollary}

\newtheorem{ex}[dfn]{Example}

\newtheorem{defn}[dfn]{Definition}

\def\Lra{\Leftrightarrow}
\def\no{\noindent}

\def\0{\emptyset}
\def\half{\frac{1}{2}}
\def\proof{\par\medskip\noindent{\it Proof: }}

\def\>{\rangle}
\def\<{\langle}

\def\C{\mathbb C}
\def\t{\widetilde}
\def\R{\mathbb R}

\def\Z{\mathbb Z}
\def\E{{\mathbb E}}
\def\p{{\mathfrak p}}

\def\h{\mathbb H}
\def\s{\mathbb S}

\def\eps{\epsilon}

\def\Ups{\Upsilon}

\def\su{{\mathfrak{su}}}

\def\ha{\widehat}
\def\fus{\circledast}

\def\th{\theta}
\def\al{\alpha}

\def\ga{\gamma}
\def\Ga{\Gamma}

\def\Si{\Sigma}

\def\la{\lambda}

\def\8{\infty}

\def\om{\omega}
\def\D{\partial}

\def\k{{\mathfrak k}}

\def\g{{\mathfrak g}}

\def\1{\hbox{\bf 1}}

\def\Hom{\mbox{\rm Hom}}
\begin{document}
\title{The Symplectic
       Geometry of Polygons in the 3-sphere}
\author[Thomas Treloar]{Thomas Treloar$^\dag$}
\address{Department of Mathematics\\
     University of Maryland\\
     College Park, MD 20742, USA}
\email{txt@math.umd.edu}
\thanks{$^\dag$\ Research partially supported by NSF grant DMS-98-03518.}
\date{\today}

\begin{abstract}
We study the symplectic geometry of the moduli spaces
$M_r=M_r(\s^3)$ of closed n-gons with fixed side-lengths in the
3-sphere.  We prove that these moduli spaces have symplectic
structures obtained by reduction of the fusion product of $n$
conjugacy classes in $SU(2)$, denoted $C_r^n$, by the diagonal
conjugation action of $SU(2)$. Here $C_r^n$ is a quasi-Hamiltonian
$SU(2)$-space.  An integrable Hamiltonian system is constructed on
$M_r$ in which the Hamiltonian flows are given by bending polygons
along a maximal collection of nonintersecting diagonals. Finally,
we show the symplectic structure on $M_r$ relates to the
symplectic structure obtained from gauge-theoretic description of
$M_r$. The results of this paper are analogues for the 3-sphere of
results obtained for $M_r(\h^3)$, the moduli space of n-gons with fixed
side-lengths in hyperbolic 3-space \cite{KMT}, and for $M_r(\E^3)$,
the moduli space of n-gons with fixed side-lengths in $\E^3$
\cite{KM1}.
\end{abstract}

\maketitle

\section{Introduction}
In this paper we study the symplectic geometry of the space of
polygons in $\s^3$ with fixed side-lengths modulo the group of
isometries.  We denote this moduli space by $M_r=M_r(\s^3)$. This
paper is continuation of \cite{KM1} and \cite{KMT}, which studied
the polygonal linkages in Euclidean 3-space and hyperbolic
3-space, respectively.

An (open) $n$-gon $P$ in $\s^3$ is an ordered $(n+1)$-tuple
$(x_1,...,x_{n+1})$ of points in $\s^3 \subset \C^2$ called the
vertices.  We join the vertex $x_i$ to the vertex $x_{i+1}$ by the
unique geodesic segment $e_i$, called the $i$-th edge (here we
must make the restriction $x_i$ and $x_{i+1}$ are not antipodal
points).  We let $Pol_n$ denote the space of $n$-gons in $\s^3$.
An $n$-gon is said to be closed if $x_{n+1} = x_1$.  We let
$CPol_n$ denote the space of closed $n$-gons. The group $G= SU(2)
\times SU(2)$ acting on $\s^3$ by $g \cdot x = g_1 x g_2^{-1}, \,
x\in \s^3$, $g = (g_1, g_2) \in G$, is the group of isometries of
$\s^3$.  Two $n$-gons $P=(x_1,...,x_{n+1})$ and
$P'=(x'_1,...,x'_{n+1})$ are equivalent if there exists $g \in G$
such that $g\cdot P = P'$, that is $g \cdot x_i = x'_i$, for all
$1\leq i \leq n+1$.

Let $r=(r_1,...,r_n) \in \R_+^n$ be an $n$-tuple of positive
numbers with $r_i < \pi$ for $1 \leq i \leq n$.  We denote by
$\t{N}_r$ the space of open $n$-gons in which the side $e_i$ a has
fixed length $d(x_i,x_{i+1}) = r_i$. We then let $\t{M}_r =
\t{N}_r \cap CPol_n, \, N_r = \t{N}_r/G,$ and $M_r = \t{M}_r/G$.
This paper examines the symplectic geometry of the space $M_r$.

We have $G = SU(2) \times SU(2)$, $K$ is the diagonal subgroup in
$G$, and $P= G/K$ which we identify with $SU(2)$.  We equip
$G,K,P$ with the quasi-Poisson structures associated to the
standard Manin pair $(\g, \k)$, where $\g = \{ (x,y) \in
\mathfrak{su}(2) \oplus \mathfrak{su}(2)\}$ and $\k = \{ (x,x)\in
\g : x \in \mathfrak{su}(2)\}$.

The main theorem of this paper is:

\begin{thm}
The space $M_r$ is a symplectic manifold with the symplectic
structure obtained from reduction of the fusion product of $n$
conjugacy classes in $SU(2)$, $C_{r_1} \fus \cdots \fus C_{r_n}$,
by the diagonal dressing action (conjugation) of the quasi-Poisson
Lie group $K$.
\end{thm}

\medskip

We are also interested in finding an integrable system on $M_r$.
We denote by $d_{ij}$ a geodesic connecting the vertices
$x_i$ and $x_j$ (we always assume $i<j$), which we call a
diagonal. Let $\ell_{ij}$ be the length of the diagonal $d_{ij}$.
Then $\ell_{ij}$ is a continuous function on $M_r$, but it is not
smooth when either $\ell_{ij}= 0$ or $\ell_{ij}= \pi$.  If
$d_{ij}$ and $d_{km}$ are nonintersecting diagonals, then
$$\{\ell_{ij}, \ell_{km} \} = 0.$$ \no By considering a maximal
collection of nonintersecting diagonals, we obtain $\half
dim(M_r)$ Poisson commuting Hamiltonians.

The Hamiltonian flow $\Psi^t_{ij}$ associated to a $\ell_{ij}$ has
the following nice description. Separate the polygon into two
pieces via the diagonal $d_{ij}$, the Hamiltonian flow is given by
leaving one piece fixed while rotating the other piece about the
diagonal at constant angular velocity 1.  The flow $\Psi^t_{ij}$
is called the ``bending flow'' along the diagonal $d_{ij}$.

\medskip

\no The paper is organized ad follows:

In section 2, we give background material for Manin pairs and
quasi-Poisson Lie groups.

In section 3, we define a symplectic structure on $M_r$ by
quasi-Hamiltonian reduction on the fusion product of conjugacy
classes.

In section 4, we study the Hamiltonians $\ell_{ij}$ and their
associated Hamiltonian flows.

In section 5, we study the an action of the pure braid group on
$M_r$ given by the time 1 Hamiltonian flows of a certain family of
functions.

In section 6, we relate the symplectic form on $M_r$ to symplectic
form given on the relative character varieties on $n$-punctured
2-spheres.

\medskip

We note that the moduli spaces of polygons in the spaces of
constant curvature give examples of completely integrable systems
obtained from the theory of Manin pairs associated to a compact
simple Lie group \cite{AMM2}.  The Manin pairs corresponding to
the various moduli spaces are:
\begin{itemize}
\item $\big(\su(2) \ltimes \su(2)^*, \su(2)\big)$ for polygons in the zero curvature
space (Lie-Poisson theory);
\item $\big(\mathfrak{sl}_2(\C) = \su(2)^{\C}, \su(2)\big)$ for polygons in negative curvature space (Poisson-Lie theory);
\item $\big(\su(2) \oplus \su(2), \su(2)\big)$ for polygons in positive curvature space (quasi-Poisson Lie theory).
\end{itemize}

\begin{center}
\textbf{Acknowledgments}
\end{center}

\normalsize The author would like to thank John Millson for
introducing him to the symplectic geometry of polygons and for
numerous fruitful discussions.  Thanks are due to Bill
Goldman for many useful conversations.  The author would also like to
thank Eckhard Meinrenken bringing to his attention \cite{AKS} and
Propostion \ref{fusion}.

\section{Manin Pairs and quasi-Poisson Lie groups}

\subsection{quasi-Poisson Structures}
In this section, we let $K$ be any compact simple Lie group with
Lie algebra denoted by $\k$.  Let $G= K \times K$ be the double of
$K$ with Lie algebra $\g = \k \oplus \k$.  The Killing form on
$\k$, which we denote by $(,)$, defines a nondegenerate bilinear
form $B(,)$ on $\g$ given by $$ B((X_1,X_2), (Y_1,Y_2)) = (X_1,
Y_1) - (X_2,Y_2), \, \textnormal{for}\, (X_1,X_2), (Y_1,Y_2) \in
\g. $$

If we now let $K$ denote the diagonal subgroup of $G$ then
its Lie algebra $\k$ is a maximal isotropic subalgebra of $\g$.
The pair $(\g,\k)$ is a Manin pair.  We will construct a
quasi-Poisson Lie group structure on $G$ associated to the Manin
pair $(\g,\k)$ which restricts to a (trivial) quasi-Poisson Lie
group structure on $K$.  For background on quasi-Poisson Lie
groups, quasi-Poisson structures, Manin pairs, etc. we refer the
reader to \cite{AKS}, \cite{Le}, \cite{KS1}, \cite{KS2}.

Let $\p = \{(\half X,-\half X) \in \g\}$ be the anti-diagonal in
$\g$.  Then $\p$ is an isotropic complement of $\k$.  Note that
$\p$ is not a Lie subalgebra of $\g$ ($[\p, \p]\subset \k$), so
the triple $(\g,\k,\p)$ is a Manin quasi-triple, rather than a
Manin triple which arises in the theory of Poisson Lie groups.  We
call this triple $(\g, \k, \p)$ the standard Manin quasi-triple.

A Manin quasi-triple gives rise to a Lie quasi-bialgebra $(\k , F ,
\varphi)$.  We can identify $\p$ with $\k^*$ via the bilinear form of
$\g$.  The cobracket on $\k$ is a map $F:\k \to \k \wedge \k$ which is
the transpose of the map from $\p \wedge \p \to \p$, also denoted by
$F$, defined by
$$
F(\xi,\eta) = \rho_\p [ \xi, \eta], \, \xi,\eta \in \p.
$$
We can also define the element $\varphi \in \wedge^3 \k$ by the map $\p
\wedge \p \to \k$ given by
$$
\varphi(\xi,\eta) = \rho_\k [ \xi, \eta], \, \xi,\eta \in \p.
$$
For the Manin quasi triple $(\g , \k ,\p)$ given above, we have $F =
0$ and $\varphi = \frac{1}{24} \sum_{ijk} f^i_{jk} e_i \wedge e_j
\wedge e_k$, where $[e_j , e_k] = \sum_i f^i_{jk} e_i$.

We can also identify $\g$ with $\k \oplus \k^*$ via the bilinear
form B(,).  The canonical $r$-matrix on
$\g$ associated to the Manin quasi-triple $(\g,\k,\p)$ is an element
$r_\g \in \g \otimes \g$ defined by the map
$r_\g :\g^* \to \g$ given by $r_\g(\xi,X)=(0,\xi)$ where $X \in \g$ and
$\xi \in \g^*$.  Let $\{e_i\}$ be an orthonormal basis of $\k$ and $\{\varepsilon^i\}$ be the dual basis in
$\k^*$, then
$$
r_\g = \sum_{i} e_i \otimes \varepsilon^i.
$$

The multiplicative 2-tensor $w_G = dL_g r_\g - dR_g r_\g$ actually defines a
bivector on $G$, since the symmetric part of $r_\g$ is a multiple of
the bilinear form $B(,)$ on $\g$. $w_\g$ gives us a quasi-Poisson Lie
group structure on $G$.  $w_\g$ naturally restricts to the trivial
bivector on the subgroup $K \subset G$.  There is also a natural
projection of $w_\g$ to $G/K = P$, which can identified with $K$, via
the map $p:G \to P$ defined by
$p(g_1,g_2)= g_1g_2^{-1}$. The bivector $w_P$ is given by

$$
w_P= \half \sum_{i} e_i^\la \wedge e_i^\rho.
$$
Here $e_i^\la$ ($e_i^\rho$) denotes the left-invariant
(resp. right-invariant) vector field on $P$ with value $e_i$ at the
identity.  We will use this notation for vector fields on $P$
throughout the rest of the paper.
Note that $w_P$ is not multiplicative, so $P$ is not a quasi-Poisson
Lie group.  We will see that in the next section that $P$ is the
target space of a generalized moment map.

\subsection{Moment map and reduction}
The action of $G$ on itself is by left multiplication induces an
action of K on $P$, the dressing action,  which is given by
conjugation.

We denote by $x_{M}$ the vector field, more generally the multivector
field,  on $M$ induced by the action of $K$ on $M$ and  $x \in \k$
satisfying
$$
(x_{M} f)(m) = \frac{d}{dt}|_{t=0} f(\exp (-tx) \cdot m)
$$
where $f \in C^\8(M)$ and $m \in M$.  This is a Lie algebra
homomorphism, i.e. $[x_M, y_M] = [x,y]_M$ for $x,y \in \k$.

We have the following definition of a quasi-Poisson action.

\begin{defn} \label{2.1}
Let $(K, w_K, \varphi)$ be a connected quasi-Poisson Lie group acting on
a manifold $M$ with bivector $w_M$.  The action of $K$ on $M$ is said
to be a quasi-Poisson action if and only if
\begin{itemize}
\item[(i)] $\half [w_M, w_M] = \varphi_M$
\item[(ii)] $\mathcal{L}_{x_M} w_M = -(F(x)_M)$
\end{itemize}
for all $x \in \k$.
\end{defn}

The dressing action of $K$ on $P$ is a quasi-Poisson action.  There is
also a notion of a generalized moment map associated to a
quasi-Poisson action.

\begin{defn}
A map $\mu:M \to P$, equivariant with respect to the action of
$K$ on $M$ and the dressing action of $K$ on $P$, is called a moment
map for the action of $K$ on $(M,w_M)$ if, on any open subset of $M$,
$$
w^\sharp(\mu^* \al_x) = x_M.
$$
Here $\al_x \in \Omega^1(P)$ is defined by $<\al_x , \xi_P> = -(x,
\xi)$ for $x \in \k$ and $\xi \in \p$.
\end{defn}

\begin{defn}
The action of $K$ on $M$ is called quasi-Hamiltonian if it admits
a moment map. A quasi-Hamiltonian space is a manifold with
bivector on which a quasi-Poisson Lie group acts by a
quasi-Hamiltonian action.
\end{defn}

The following lemma will be useful in this paper for the proofs of
Proposition \ref{fusion} and Theorem \ref{symplec}.

\begin{lem} \label{quasiham}
Let $(M,w_M)$ be a manifold with bivector on which the compact
simple Lie group $K$ act in a quasi-Poisson manner.  Then
$(M,w_M)$ is a quasi-Hamiltonian space if and only if there exists
a map $\mu:M \to P$ which is equivariant with respect to action of
$K$ on $M$ and the action of $K$ on $P$ by conjugation which
satisfies $$ w^{\sharp} (\mu^*  (x,\theta)) = \half ( (1_\k +
Ad_\mu)x)_M $$ for all $x \in \k$. Here $w^\sharp : T^* M \to
T_*M$ is given by $w^\sharp(\al) = w(\al, \cdot)$ for $\al \in
T^*M$, and $\theta:T_*K \to \k$ is the left-invariant
Maurer-Cartan on $K$.  For $K$ a matrix group $\theta = k^{-1}
dk$.
\end{lem}
\proof
See \cite[Proposition 5.33]{AKS}.
\qed
\begin{ex}
The basic example of a quasi-Hamiltonian space is the space $P$.
The action of $K$ on $P$ is the dressing action and the associated
moment map is the identity map.  The bivector on $P$ is given by
$w_P= \half \sum_{i} e_i^\la \wedge e_i^\rho$.
\end{ex}

In general, any $K$-invariant embedded submanifold of $P$ is also
a quasi-Hamiltonian space with moment map given inclusion.

\begin{ex} \label{conjclass}
Let $(\g, \k, \p)$ be the standard Manin quasi-triple.  Let $C
\subset P$ be a conjugacy class in $P$. The action of $K$ on $C$
given by conjugation is a  quasi-Poisson action. The momentum map
associated to this action of is the inclusion map (i.e. $\mu:C \to
P$ given by $\mu(g) = g$). Since the bivector $w_P$ is
$K$-invariant, the bivector on $C$ is given by the restriction
$w_P|_C$
\end{ex}

Even though a quasi-Hamiltonian space $(M,\mu, w_M)$ is not in
general a Poisson manifold, $\half [w_M, w_M] = \varphi_M$, there
is still a notion of reduction to a symplectic manifold.

\begin{lem} \label{symplec}
Let $(M, w_M ,\mu)$ be a quasi-Hamiltonian space such that the
bivector $w_M$ is everywhere nondegenerate.  Assume $M/G$ is a
smooth manifold in a neighborhood $U$ of $p(x_0)$, where $p:M \to
M/G$ and $x_0 \in M$. Let $x \in M$ be such that $p(x) \in U$ and
$s=\mu(x) \in D/G$ is a regular value of the moment map $\mu$.
Then the symplectic leaf through $p(x)$ in the Poisson manifold
$U$ is the connected component of the intersection with $U$ on the
projection of the manifold $\mu^{-1}(s)$.
\end{lem}

\proof
See \cite[Theorem 5.5.5]{AKS}

\subsection{Fusion product of quasi-Poisson manifolds}

Given quasi-Hamiltonian spaces $M_1$ and $M_2$ each acted on by
$K$ with associated moment maps $\mu_1:M_1 \to P$ and $\mu_2:M_2
\to P$, it is not true that $M_1 \times M_2$ with the product
bivector structure is a quasi-Hamiltonian K-space with the action
being the diagonal action of $K$ on $M_1 \times M_2$. We can
define a new bivector on $M_1 \times M_2$ such that diagonal
action is a quasi-Poisson action with respect to this new
bivector.  $M_1 \times M_2$ with this bivector is called the
fusion product and is due to \cite{AKSM}.

As defined in the previous section, the subscript $M$ denotes the vector field, or
multivector field, induced by the action of $K$ on $M$.

\begin{prop} \label{fusion}
Let $(M_1, w_1, \mu_1)$ and $(M_2, w_2, \mu_2)$ be
quasi-Hamiltonian K-spaces in the sense of \cite{AKS}. Then $M =
M_1 \times M_2$ with the action of $K$ on $M$ given by the
diagonal action, bivector on $M$ given by $$ w_M = w_1 + w_2 +
\half \sum_j (e_j)_{M_1} \wedge (e_j)_{M_2} $$ and moment map $\mu
= \mu_1\mu_2$ is a quasi-Hamiltonian $K$-space.  Recall $\{e_i\}$
is an orthonormal basis of $\k$.   $M$ with this structure is
called the fusion product of $M_1$ and $ M_2$ and is denoted by $M
= M_1 \circledast M_2$.
\end{prop}

\proof
We begin by showing the diagonal action of $K$ on $(M,w_M)$ is a
quasi-Poisson action.  For this we need to show,

\begin{itemize}
\item[(i)] $\half [w_M, w_M] = \varphi_M$
\item[(ii)] $\mathcal{L}_{x_M} w_M = 0$.
\end{itemize}

\medskip

We will then show that $\mu:M_1 \times M_2 \to P$ given above is the
moment map associated to the diagonal action.

It is a straightforward calculation to show $(i)$:

\begin{eqnarray*}
\half \Big[w_M, w_M\Big] & = & \half \Big[w_1 + w_2 + \half \sum_j
(e_j)_{M_1} \wedge (e_j)_{M_2}, w_1 + w_2 + \half \sum_k (e_k)_{M_1} \wedge
(e_k)_{M_2}\Big] \\
    & = & \half \Big[w_1, w_1\Big] + \half\Big[w_2, w_2\Big] + \Big[w_1 + w_2, \half
\sum_{j=1}^n (e_j)_{M_1} \wedge (e_j)_{M_2}\Big] \\
    &   & \, + \half\Big[ \half \sum_j
(e_j)_{M_1} \wedge (e_j)_{M_2}, \half \sum_k (e_k)_{M_1} \wedge
(e_k)_{M_2}\Big]\\
    & = & \half \Big[w_1, w_1\Big] + \half\Big[w_2, w_2\Big] +
\Big[w_1 +w_2 , \sum_j (e_j)_{M_1} \wedge (e_j)_{M_2}\Big] \\
    &   & + \frac{1}{8}
 \sum_{j,k} \Big(\Big[(e_j)_{M_1},(e_k)_{M_1}\Big] \wedge (e_j)_{M_2} \wedge
(e_k)_{M_2} + \Big[(e_j)_{M_2},(e_k)_{M_2}\Big] \wedge (e_j)_{M_1}
\wedge (e_k)_{M_1} \Big)\\
\end{eqnarray*}

But $ \half \Big[w_i, w_i\Big] =\varphi_{M_i}$ for $i= 1,2$ since the
    $K$- actions on $M_1$ and $M_2$ are quasi-Poisson actions.  Also, we
    have $[(e_k)_{M_i}, w_i] = \mathcal{L}_{(e_k)_{M_i}} w_1 =
    -\Big(F(e_k)\Big)_{M_i}$ where $F: \k \to \wedge^2 \k$ is the
    cobracket.  But $F\equiv 0$ for the standard quasi-Poisson Lie
    group $K$ we have,  thus $[(e_k)_{M_i}, w_i] = 0$.  Let $f^i_{jk}$
    denote the structure constants on $\k$. The above equations
    then become

\begin{eqnarray*}
    & = & \varphi_{M_1} + \varphi_{M_2} + 0 + \frac{1}{8}
\sum_{j,k} \Big[e_j,e_k\Big]_{M_1} \wedge (e_j)_{M_2} \wedge
(e_k)_{M_2} \\
    &   & \, + \frac{1}{8} \sum_{j,k} \Big[e_j,e_k\Big]_{M_2} \wedge
(e_j)_{M_1} \wedge (e_k)_{M_1}\\
    & = & \frac{1}{24} \sum_{ijk} f_{jk}^i (e_i)_{M_1} \wedge (e_j)_{M_1}
\wedge (e_k)_{M_1} + \frac{1}{24} \sum_{ijk} f_{jk}^i (e_i)_{M_2}
\wedge (e_j)_{M_2} \wedge (e_k)_{M_2} \\
    &   & \, + \frac{1}{8} \sum_{ijk} f_{jk}^i (e_i)_{M_1}
\wedge (e_j)_{M_2} \wedge (e_k)_{M_2} + \frac{1}{8} \sum_{ijk}
f_{jk}^i (e_i)_{M_2} \wedge (e_j)_{M_1} \wedge (e_k)_{M_1} \\
    & = & \frac{1}{24} \sum_{ijk} f_{jk}^i \Big((e_i)_{M_1}+(e_i)_{M_2}\Big)
\wedge \Big((e_j)_{M_1} + (e_j)_{M_2}\Big) \wedge \Big((e_k)_{M_1} + (e_k)_{M_2}\Big)\\
    & = & \frac{1}{24} \sum_{ijk} f_{jk}^i (e_i)_M \wedge (e_j)_M
\wedge (e_k)_M\\
    & = & \varphi_M
\end{eqnarray*}

To show $(ii)$, we again use $\mathcal{L}_{(e_k)_{M_i}} w_{M_i} = 0$.

\begin{eqnarray*}
\mathcal{L}_{(e_k)_M} w_M & = & \mathcal{L}_{(e_k)_{M_1} + (e_k)_{M_2}}  \Big(w_1 + w_2
    + \sum(e_j)_{M_1} \wedge (e_j)_{M_2}\Big) \\
    & = & \mathcal{L}_{(e_k)_{M_1} + (e_k)_{M_2}} \Big( \sum(e_j)_{M_2} \wedge
    (e_j)_{M_2}\Big) \\
    & = & \sum \Big[(e_k)_{M_1}, (e_j)_{M_1}\Big] \wedge
    (e_j)_{M_2} - \sum \Big[(e_k)_{M_2}, (e_j)_{M_2}\Big] \wedge
    (e_j)_{M_1}\\
    & = & \sum_{i,j} C_{kj}^i (e_i)_{M_1} \wedge (e_j)_{M_2} -
\sum_{i,j} C_{kj}^i (e_i)_{M_2} \wedge (e_j)_{M_1}\\
    & = & 0
\end{eqnarray*}

We next use Lemma \ref{quasiham} to show that $\mu = \mu_1 \mu_2: M_1
\times M_2 \to P$ is indeed the moment map associated to the diagonal action.

\begin{eqnarray*}
w^{\sharp} (\mu^*(x,\theta)) & = & w^{\sharp} ((\mu_1\mu_2)^*
(x,\theta))\\
    & = & w^{\sharp} ((x,\mu_2^*\theta + Ad_{\mu_2^{-1}} \mu_1^*
\theta))\\
    & = & w^{\sharp} (\mu_2^*(x, \theta) + \mu_1^* (Ad_{\mu_2}x,
\theta))\\
    & = & w_1^\sharp\Big(\mu_1^*(Ad_{\mu_2}x, \theta)\Big) +
w_2^\sharp\Big(\mu_2^*(x, \theta)\Big) + \half \sum_j \Big(
(\mu_1^*(Ad_{\mu_2}x,\theta))(e_j)_{M_1}\Big)\, (e_j)_{M_2} \\
    &   & \, - \half \sum_j \Big
( (\mu_2^*(x,\theta))(e_j)_{M_2}\Big)\, (e_j)_{M_1}\\
\end{eqnarray*}

\no $(M_i, w_i)$ is a quasi-Hamiltonian space with moment map
$\mu_i: M_i \to P_i$, so we have by Lemma \ref{quasiham}

$$
 w_i^\sharp\Big(\mu_i^*(x, \theta)\Big) = \half((1+
 Ad_{\mu_i})x)_{M_i}.
$$

\no We can also see that

\begin{eqnarray*}
\sum_i \Big( (\mu_j^*(x,\theta))(e_i)_{M_j}\Big)\, (e_i)_{M_k} &
= & \sum_i (x, Ad_{\mu_j^{-1}} e_i - e_i)(e_i)_{M_k}\\
   & = & \sum_i (Ad_{\mu_j}x - x, e_i)(e_i)_{M_k} \\
   & = & (Ad_{\mu_j}x - x)_{M_k}
\end{eqnarray*}

\no So the above becomes

\begin{eqnarray*}
w^{\sharp} (\mu^*(X,\theta))  & = & \half(Ad_{\mu_2} + Ad_{\mu_1\mu_2}X)_{M_1} + \half(1 +
Ad_{\mu_2}X)_{M_2} + \half (Ad_{\mu_1\mu_2}X -
Ad_{\mu_2}X)_{M_2}\\
    &   & \, -\half (Ad_{\mu_2}X - X)_{M_1}\\
    & = & \half((1+ Ad_{\mu_1\mu_2})X)_{M_1} + \half((1+
Ad_{\mu_1\mu_2})X)_{M_2}\\
    & = & \half((1+ Ad_{\mu_1\mu_2})X)_{M}
\end{eqnarray*}

\qed

\begin{rem}
It is a quick calculation to show the fusion product is associative,
that is $M_1 \circledast (M_2 \circledast M_3) \simeq  ( M_1 \circledast
M_2) \circledast M_3$.  The bivector is given by
$$
w = w_1 + w_2 + w_3
+ \half \sum_i (e_i)_{M_1} \wedge (e_i)_{M_2} + \half \sum_i
(e_i)_{M_1} \wedge (e_i)_{M_3} + \half \sum_i (e_i)_{M_2} \wedge
(e_i)_{M_3}.
$$
\end{rem}

The quasi-Hamiltonian space we are most interested in for this
paper is the fusion product of $n$ conjugacy classes in $P$.
Recall from Example \ref{conjclass} that $C_{r_i} \subset P$ is a
quasi-Hamiltonian space with action given by conjugation and the
associated moment map given by inclusion.  The fusion product of
$n$ conjugacy classes $C^n_r = C_{r_1} \fus \cdots \fus C_{r_n}$,
$r=(r_1, ..., r_n) \in \R_+$ is also a quasi-Hamiltonian space
with action given by the diagonal conjugation and moment map
$\t{\mu}: M \to P$ given by multiplication,  $\t{\mu}(g_1, g_2,
..., g_n) = g_1 g_2 \cdots g_n$. The bivector on this space is
given by $$ \t{w} = \half  \sum_{i=1}^n  \sum_k \Big(e_k^\la
\wedge e_k^\rho \Big)_i + \half \sum_{i<j}^n  \sum_k \Big(e_k^\la
- e_k^\rho \Big)_i \wedge \Big(e_k^\la - e_k^\rho \Big)_j $$ where
the subscripts $i,j$ denote the vector field on $C_{r_i}, C_{r_j}
\subset C^n_r$.

\subsection{Poisson bracket on $C^\8(P^n)^K$}

For a general quasi-Hamiltonian space $(M,w_M)$, the bracket on
$C^\8(M)$ defined by the bivector $w_M$ is not a Poisson bracket.
This is easy to see since the Shouten bracket $[w_M,w_M] =
\varphi_M$ is an invariant trivector field.  The bracket does
however define a Poisson bracket when we restrict to the  space
$C^\8(M)^K$ of smooth K-invariant functions on $M$.

\begin{lem}
Let $K$ be a connected quasi-Poisson Lie group acting on a manifold
$(M,w_M)$ in a quasi-Poisson manner.  Then the bivector $w_M$ defines
a Poisson bracket on the space $C^\8(M)^K$ of the smooth $K$-invariant
functions in $M$.
\end{lem}

\proof
See \cite[Theorem 4.2.2]{AKS}
\qed

\medskip

For  $\psi\in C^\8 (P^n)$ we define
$$
D_i \psi : P^n \to \k_i , \quad D_i' \psi : P^n \to \k_i
$$
as follows. Let $g=(g_1,...,g_n)\in P^n$ and $x=(x_1, ..., x_n) \in \k^n$, then
$$
d_i \psi_g(x^\rho) = (D_i \psi,x)  = \frac{d}{dt}|_{t=0}
\psi(g_1,...,e^{t x_i}g_i, ..., g_n)
$$
$$
d_i \psi_g(x^\la) = (D_i '\psi,x)  =\frac{d}{dt}|_{t=0}
\psi(g_1,..., g_i e^{t x_i}, ..., g_n).
$$

Here $(,)$ is the Killing form extended to $\k^n$ by $(x,y) =
\sum_{i=1}^n(x_i, y_i)$ for $x, y \in \k^n$.

\begin{rem}
It is easy to see that
$$
Ad_{g_i} D_i' \psi(g) = D_i\psi
$$
\end{rem}

We also define
$$\Psi_j(g) = \sum_{i=1}^{j-1} \Big[D_i \psi(g) - D'_i \psi(g)\Big]
+ D_j \psi(g)
$$

We now define the  Poisson bracket on $C^\8(P^n)^K$.
\begin{prop} \label{bracket}
Let $\phi,\psi \in C^\8(P^n)^K$ then
$$
\{\phi,\psi \}(g) = \sum_{j=1}^n \Big(D'_j\phi(g) -
D_j\phi(g), \Psi_j(g)\Big)
$$
\end{prop}
\proof

Let us first note that for $x,y \in \k$  $\sum_i (x,e_i)(y,e_i)
= (x,y)$.  Now,

\begin{eqnarray*}
\{\varphi,\psi \}(g) &  = & w(d\varphi,d\psi) \\
    & = & \half \sum_{i=1}^n \sum_k \Big(e_k^\la \wedge e_k^\rho
\Big)_i(d\phi,d\psi) + \half \sum_{i<j}^n \sum_k \Big((e_k^\la -
    e_k^\rho)_i
\wedge (e_k^\la - e_k^\rho)_j \Big)(d\phi,d\psi) \\
    & = &  \half \sum_{i=1}^n \sum_k d_i\phi(e_k^\la) d_i \psi
    ( e_k^\rho) - d_i\phi(e_k^\rho) d_i \psi( e_k^\la) \\
    &   & \: + \half \sum_{i<j}^n \sum_k d_i \phi(e_k^\la - e_k^\rho )
d_j \psi (e_k^\la - e_k^\rho ) - d_j \phi(e_k^\la - e_k^\rho ) d_i
\psi (e_k^\la - e_k^\rho ) \\
    & = & \half \sum_{i=1}^n \sum_k \Big(D'_i\phi,e_k \Big)
    \Big(D_i \psi, e_k\Big) - \Big(D_i\phi,e_k \Big) \Big(D'_i
    \psi, e_k\Big)  \\
    &   & \: + \half \sum_{i<j} \sum_k \Big(D'_i \phi- D_i\phi,e_k
    \Big)  \Big(D'_j \psi- D_j\psi,e_k \Big) - \Big(D'_j \phi-
    D_j\phi,e_k \Big)  \Big(D'_i \psi- D_i\psi,e_k \Big) \\
    & = & \half \sum_{i=1}^n \Big(D'_i\phi, D_i \psi \Big) -
    \Big(D_i\phi,D'_i \psi \Big)  \\
    &   & \: + \half \sum_{i<j} \Big(D'_i \phi - D_i\phi,
    D'_j \psi -  D_j\psi \Big) - \Big(D'_j \phi-
    D_j\phi, D'_i \psi - D_i\psi \Big) \\
    & = & \half \sum_{i=1}^n \Big(D'_i\phi, D_i \psi \Big) -
    \Big(D_i\phi,D'_i \psi \Big)  \\
    &   & \: + \half \sum_{i<j} \Big(D'_i \phi - D_i\phi,
    D'_j \psi -  D_j\psi \Big) - \sum_{i>j} \Big(D'_i \phi-
    D_i\phi, D'_j \psi - D_j\psi \Big) \\
\end{eqnarray*}

But since $\psi \in C^\8(P^n)^K$ is $K$-invariant, a quick calculation
shows

$$
\sum_{i=1}^n [D_i \psi - D'_i \psi] = 0
$$

Using this fact and also that $(D'_i\phi, D'_i\psi) = (D_i\phi,
D_i\psi)$ for all $i$, we can rewrite the above as,
\begin{eqnarray*}
\{\phi,\psi\} & = & \half \sum_{i=1}^n \Big(D_i' ,\phi - D_i\phi, D_i \psi
    + D'_i \psi \Big)  \\
    &    & \: - \half \sum_{i\geq j} \Big(D'_i \phi - D_i\phi,
    D'_j \psi -  D_j\psi \Big) - \half \sum_{i > j} \Big(D'_i \phi-
    D_i\phi, D'_j \psi - D_j\psi \Big) \\
    & = & \sum_{i=1}^n \Big(D'_i \varphi - D_i\varphi, \Psi_i \Big) \\
\end{eqnarray*}
\qed

From the above Proposition we can also define the Hamiltonian
vector field $X_\psi$ associated to $\psi \in C^\8(P^n)^K$ by
$X_\psi = w^\sharp(d\psi)$.

\begin{cor} \label{hvectC_r}
The Hamiltonian vector field $X_\psi(g)=((X_1(g), ..., X_n(g))$
associated to the $K$-invariant function $\psi \in C^\8(P^n)^K$ is
given by $$ X_j(g) = dL_{g_j} \Psi_j - dR_{g_j} \Psi_j, \, 1 \leq
j \leq n. $$ and $ g=( g_1, g_2, ... , g_n)$.
\end{cor}
\proof
We use the convention $\{\phi,\psi\} = d\phi(X_\psi) = \sum_{j=1}^n
d_j \varphi((X_j(g))$.  Proposition \ref{bracket} gives us
\begin{eqnarray*}
d\phi(X_\psi(g)) & = & \{\phi, \psi\} \\
    & = & \sum_{j=1}^n\Big(D'_j\phi - D_j \phi, \Psi_j
    \Big)\\
    & = &  \sum_{j=1}^n d_j \phi (dL_{g_j}\Psi_j) - d_j \phi
    (dR_{g_j}\Psi_j)\\
    & = &  \sum_{j=1}^n d_j \phi (dL_{g_j}\Psi_j - dR_{g_j}\Psi_j)
\end{eqnarray*}
\qed

\section{The symplectic structure on $M_r(\s^3)$}
Throughout the rest of the paper, we let $G
= SU(2) \times SU(2)$, $K = SU(2)$, and $P \simeq SU(2)$.  In this
section, we will define a symplectic structure on $M_r$ obtained from
the reduction of the fusion product of conjugacy classes to a
symplectic manifold.

Recall, we defined $Pol_n(*)$ to be the open $n$-gons in $\s^3$ with
side-length less than $\pi$, so that we can choose an unique geodesic
between vertices.
The map $\Phi : P^n \to Pol_n(*) \subset (S^3)^n$ defined by
$$
\Phi(g) = ( *, g_1 *, g_1 g_2 *, ..., g_1 g_2 \cdots g_n *)
$$
is a diffeomorphism.

\begin{prop}
The map $\Phi$ is a $K$-equivariant diffeomorphism  where $K$ acts on
$P^n$ by the dressing action (diagonal conjugation) and on $Pol_n(*)$ by the
diagonal action on $(\s^3)^n$.
\end{prop}
\proof
$* \in P$ is an element in $P$ which is fixed by the
$K$-action, that is $Ad_k(*) = *$ for all $k \in K$.  For $k\in K$ and
$g\in P^n,\, k\cdot p=(Ad_k g_1, ..., Ad_k g_n)$, so
\begin{eqnarray*}
\Phi(k\cdot g) & = & (*, Ad_k(g_1) *, ..., Ad_k(g_1\cdots g_n)*)\\
    & = & (Ad_k*, Ad_k(g_1 *, \cdots, Ad_k(g_1\cdots g_n*)) \\
    & = &  k \cdot (*, g_1 *, ..., g_1 \cdots g_n *).
\end{eqnarray*}
\qed

\begin{rem}
The map $\Phi$ induces a diffeomorphism from $\{g\in P^n: g_1 \cdots
g_n = 1 \}$ to $CPol(*)$.
\end{rem}

We have seen that the $K$-orbits in a quasi-Hamiltonian space are
quasi-Hamiltonian spaces.  In particular, a conjugacy class $C
\subset P$ is a quasi-Hamiltonian space.  Let $r \in \R^n$, with
$r=(r_1, ..., r_n)$.  Let $C_{r_i} \subset P$ denote the conjugacy
class in $P$ such that $r_i = d(*,g_i*) = \cos^{-1}\Big(-\half
trace(g_i)\Big)\in \R$ for all $g_i \in C_{r_i}$.

\begin{lem}
The map $\Phi$ induces a $K$-equivariant diffeomorphism from
$C_{r_1} \times \cdots \times C_{r_n}$ to  $\t{N}_r$, the space of
open $n$-gons with fixed side-lengths based at $*$, where $r_i =
d(g_1 \cdot g_i *, g_1 \cdot g_{i-1}*)$, for all $1 \leq i \leq n$.
\end{lem}
\proof
Follows from the fact that $k$ fixes side-lengths.
\qed

\begin{cor}
$\Phi$ induces a diffeomorphism  from the space  $\{g\in C^n_r: g_1 \cdots
g_n = 1 \}/K$  to $M_r$ the moduli space of closed $n$-gons in $\s^3$.
\end{cor}

In $\S 2.3$ we saw that the fusion product of $n$ conjugacy
classes in $P$, $(C^n_r,\t{\mu}, \t{w})$, is a quasi-Hamiltonian
space with the moment map $\t{\mu}$ given by multiplication.  So,
$\t{\mu}^{-1}(1)/K = \{g\in C^n_r: g_1 \cdots g_n = 1 \}/K$.  We
must determine when this restriction and quotient gives rise to
symplectic manifold.  Lemma \ref{symplec} tells us that
$\t{\mu}^{-1}(1)/K$ is a symplectic manifold when
\begin{itemize}
\item $\t{w}$ is everywhere nondegenerate on $C^n_r$
\item 1 is a regular value of $\t{\mu}$.
\end{itemize}

\no We use the following remark from \cite[Example 5.5.4]{AKS} to give
the nondegeneracy condition.
\begin{rem}
Let $K$ be a quasi-Poisson Lie group arising from the standard
quasi-triple and $(M, \mu, w)$ is a  quasi-Hamiltonian space. Then
$(M, \mu , w)$ is nondegenerate if and only if, for each $m \in
M$, $$ ker(w_m^\sharp) = \{\mu^*(x, \theta) : x \in ker (1 +
Ad_{\mu(m)})\}. $$ Here $x \in \k$.
\end{rem}

\no It follows that the fusion product of conjugacy classes is
nondegenerate.

\begin{lem}
1 is a regular value of $\t{\mu}$ if and only if $\k_g = \{x \in \k : x_{C_r^n}
= 0 \} = 0$ for all $g \in \t{\mu}^{-1}(1)$.
\end{lem}
\proof We refer to Lemma \ref{quasiham}. Let $x\in \k$. Then
$x \in (Im(d\t{\mu}|_g))^\perp$ $\Lra$ $(x, \t{\mu}^* \theta) = 0$
$\Lra$ $0 = \t{w}^\sharp((x, \t{\mu}^* \theta)) =
((1+Ad_{\t{\mu}(g)})x)_{C_r^n} = (2x)_{C_r^n}$.
\qed

A polygon is said to be degenerate if it can be contained in a
geodesic in $\s^3$.  It follows from the above lemma that if there does not
exist $g \in \t{\mu}^{-1}(1) \subset C_r^n$ such that $\Phi(g)$ is a
degenerate polygon, then 1 is a regular value of $\t{\mu}$.

\begin{thm}
The moduli space $M_r$ containing no degenerate
polygons has a symplectic structure which is the transport structure
from the moduli space $\mu^{-1}(1)/K$.
\end{thm}

\no In \S 6, we need a formula for the symplectic form on $M_r$ i in \S 6.

\begin{rem} \label{symform}
The symplectic form is given by  
$$
\t\om = \sum_{i=1}^n \om_i + \half \sum_{i=1}^n \sum_{j=i+1}^n
\left( Ad_{g_1 \cdots g_{i-1}} \bar\th_i \wedge_b Ad_{g_1 \cdots
g_{j-1}} \bar\th_j \right). $$
where $\om_i$ is the quasi-Hamiltonian 2-form on the conjugacy class
$C_i \subset SU(2)$, see \cite{AMM}, and $\bar\th_i$ is the
right-invariant Maurer-Cartan form on $C_i \subset SU(2)$.  We denote
by $\wedge_b$ the wedge product together with the killing form on $G$.
\end{rem}

\section{Bending Hamiltonians}

\subsection{Hamiltonian vector fields}

Recall, $K=SU(2)$ and $C^n_r= C_{r_1} \circledast \cdots
\circledast C_{r_n}$, where $C_{r_i} \subset P$ is a conjugacy
class in $P\simeq SU(2)$.  Let $(x,y) = -\half Tr(xy)$.  In this
section we will compute the Hamiltonian
vector fields $X_{f_j}$ associated to the functions $f_i \in
C^\8(C^n_r)^K$ given by
$$ f_j(g)= tr(g_1 \cdots g_j), \; 1\leq j\leq n. $$

See $\S 2.4$ for the definition of the  Poisson bracket on $C^\8(C^n_r)^K$.  We leave it to
the reader to verify the following lemma.

\begin{lem} \label{4.1}
\begin{eqnarray*}
D_{i+1} f_j(g) &  = & D'_{i} f_j(g), \,\, 1 \leq i \leq j-1 \\
D_1 f_j(g) & = & D'_{j} f_j(g)
\end{eqnarray*}
\no for all $1 \leq j \leq n$.
\end{lem}

\no We define $F_j: P \to \k$ by
$$ F_j(g) = \Big((g_1 \cdots g_j) - (g_1 \cdots g_j)^{-1}\Big). $$

\no We then have the following lemma.

\begin{lem} \label{4.2}
$F_j(g) = D_1 f_j(g)$
\end{lem}

\proof For $g\in C^n_r$ and $X \in \k$
\begin{eqnarray*}
(D_1 f_j(g) , X) & = & \frac{d}{dt} \Big|_{t=0} tr(e^{tX} g_1 g_2
    \cdots g_j) \\
    & = & tr(X g_1 g_2 \cdots g_j) \\
    & = & tr(g_1 g_2 \cdots g_j X)
\end{eqnarray*}
but since
$$
tr((g_1 g_2 \cdots g_j)^{-1}X) = tr((g_1 \cdots
g_j)^*X) = tr(X^*g_1 \cdots g_j) = -tr(g_1 \cdots g_j X)
$$
it follows that
\begin{eqnarray*}
tr(g_1 g_2 \cdots g_j X) & = & \half  tr\Big(((g_1 g_2 \cdots g_j)-(g_1
    \cdots g_j)^{-1}) X\Big)\\
    & = & \Big(- ((g_1 \cdots g_j) - (g_1 \cdots g_j)^{-1}), X\Big).
\end{eqnarray*}
Since $-\Big((g_1 \cdots g_j) - (g_1 \cdots g_j)^{-1}\Big) \in
\k$ and $(,)$ is a nondegenerate bilinear form, we have $D_1 f_j(g)
= - \Big((g_1 \cdots
g_j) - (g_1 \cdots g_j)^{-1}\Big) = - F_j(g)$. \qed

We have the following formula of the Hamiltonian vector fields
$X_{f_i}$.

\begin{thm}\label{hamvf}
The Hamiltonian vector field $X_{f_i}$ is has an $i$-th component
given by $$(X_{f_j}(g))_i = dR_{g_i} F_j(g) - dL_{g_i} F_j(g), \;
1 \leq i \leq j,$$ $$(X_{f_j}(g))_i = 0, \; j < i \leq n$$
\end{thm}

\proof Recall from Corollary \ref{hvectC_r}  that for $\psi \in C^\8(C^n_r)^K,
X_\psi(g)$ is given by
$$ (X_{\psi}(g))_i = dL_{g_i} \Psi_i (g)
- dR_{g_i} \Psi_i (g) $$
where $\Psi_i (g) = D_1\psi(g) -
D'_1\psi(g)+ D_2\psi(g)- \cdots - D_{i-1}\psi(g) +
D_i\psi(g)$.  This together with Lemma \ref{4.1} gives us

$$(X_{f_j}(g))_i = dL_{g_i} D_1 f_j(g) - dR_{g_i} D_1 f_j(g),  \;
1 \leq i \leq j$$ and $$(X_{f_j}(g))_i= 0, \; j < i \leq n.$$

But from Lemma \ref{4.2},  $- F_j(g) = D_1 f_j(g)$,  completing the proof.

\qed

\subsection{Commuting flows} \label{commflow}

In this section we will show the family of Hamiltonians
$\{f_j\}_{j=1}^n$ Poisson commute for $1 \leq j \leq n$.

\begin{prop}
$\{f_i,f_j\} \equiv 0$ for all $i,j$.
\end{prop}

\proof
Without loss of generality we may assume $i < j$, then by Proposition
\ref{bracket}
\begin{eqnarray*}
\{f_i,f_j\}(g) & = & \sum_{k =1}^j \Big(D'_k f_i(g) - D_k f_i(g),
    F_j(g)\Big) \\
    & = & - \Big(\sum_{k =1}^j (D'_k f_i(g) - D_k f_i(g)),
     F_j(g)\Big) \\
    & = & \Big( 0 , F_j(g)\Big) \\
    & = & 0
\end{eqnarray*}
Here we used  $\sum_{k=1}^i (D_k f_i - D'_k f_i) = 0$.
\qed

\subsection{Hamiltonian flow}
In this section we will calculate the Hamiltonian flow,
$\Phi^t_j$, associated to $f_j$.  Recall that the Hamiltonian flow
is the solution to the ODE
\begin{equation*}(*)
\begin{cases}
\frac{dg_i}{dt} = dR_{g_i} F_j(g) - dL_{g_i} F_j (g), \: 1 \leq i
\leq j \\
\frac{dg_i}{dt} = 0 , j < i \leq n
\end{cases}
\end{equation*}

\begin{lem}
$F_j(g)$ is invariant along solution curves of (*).
\end{lem}

\proof
To prove the lemma, it suffices to show that $\psi_j (g) = g_1 \cdots
g_j$ is invariant along solution curves of (*).
\begin{eqnarray*}
\frac{d}{dt}\psi_j (g(t)) & = & \frac{d}{dt} (g_1(t) g_2(t) \cdots
g_j(t)) \\
    & = & \frac{dg_1}{dt}(t)  g_2(t) \cdots g_j(t) +k_1(t)
    \frac{dg_2}{dt}(t)
     \cdots g_j(t) + \cdots + g_1(t) g_2(t) \cdots
     \frac{dg_j}{dt}(t)
   \\
    & = & [F_j(g(t)) g_1(t) - g_1(t) F_j (g(t))]g_2(t) \cdots g_j(t) +
    g_1(t) [F_j(g(t)) g_2(t) - g_2(t) F_j (g(t))] \cdots g_j(t) \\
    & = & g_1(t) g_2(t) \cdots [F_j(g(t)) g_j(t) - g_j(t) F_j
    (g(t))] \\
    & = & F_j(g(t)) g_1(t) \cdots g_j(t) - g_1(t) \cdots g_j(t)
    F_j(g(t)) \\
    & = & 0
\end{eqnarray*}
\qed

\begin{lem} The curve $\exp{\big(t F_j(g)\big)}$ is periodic with period
$2\pi / \sqrt{4-{f^2_j}}$.
\proof Left to reader.
\qed

\no We are now able to find the Hamiltonian flow $\Phi_{j}^t$.
\end{lem}
\begin{thm} The Hamiltonian flow, $\Phi_j^t$, associated to the
Hamiltonian $f_j$ given by $\Phi_j^t(g) = \big(\t{g}_1(t),
...,\t{g}_n(t)\big) $ where
\begin{equation*}\t{g}_i(t) =
\begin{cases}
Ad\big(\exp(tF_j(g))\big) g_i , \: 1 \leq i \leq j \\
g_i  , j < i \leq n.
\end{cases}
\end{equation*}
The flow is periodic with period $2\pi / \sqrt{4-{f^2_j}}$.
\end{thm}

The flows \{$\Phi_j^t$\} do not give rise to a torus action on $M_r$
since they do not have constant period.  We now look at the length functions
$\ell_j(g) =  \cos^{-1} (-\half f_j(g))$. Then

$$
d \ell_j = \frac{1}{\sqrt{4 - f^2_j}} df_j
$$
\no and
$$
X_{\ell_j} = \frac{1}{\sqrt{4 - f^2_j}} X_{f_j}.
$$

It is not difficult to see that the family of functions $\{\ell_j
\}_{j=2}^{n-1}$ also Poisson commute, but their Hamiltonian flows are not
everywhere defined.  If we restrict to the space $M'_r$ such $\ell_j
\ne 0$ or $\ell_j \ne \pi$ for all diagonals in $M_r$.  The
Hamiltonian flows $\{\Psi_j^t\}$ on $M'_r$ associated to $\{\ell_j\}$
are periodic with constant period $2\pi$ and constant angular
velocity 1.  These flows define a Hamiltonian $(n-3)$-torus action on
the space $M'_r$

\section{Braid action on $M_r$}
There exists an action of the pure braid group $P_n$ on the
manifold $M_r$ which preserves the symplectic structure.  In this
section,  we show that the generators of the pure braid group
arise as the time 1 Hamiltonian flows of the family of functions
$h_{ij}, 1 \leq i < j \leq n-1$ where $h_{ij} \in C^\8(M_r)^K$ is
defined by, $$ h_{ij}(g) = \half \Big(\cos^{-1}\big(-\half tr(g_i
g_j)\big)\Big)^2. $$

Let $C_{12}$ denote $C_1 \circledast C_2$, where $C_i \subset P$
is a conjugacy class. Let $w_{12}$ denote the quasi-Poisson
bivector on $C_{12}$.   We have the following proposition.

\begin{prop}
The diffeomorphism  $R: C_1 \circledast C_2 \to C_2 \fus C_1$ given by
$R(g_1, g_2) = (Ad_{g_1}g_2,
g_1)$ is a bivector map taking $w_{12}$ to $w_{21}$.
\end{prop}

\begin{rem}
The diffeomorphism $R':C_1 \fus C_2 \to C_2 \fus C_1$ given by $R'(g_1,g_2) =
(g_2, Ad_{g_2^{-1}} g_1)$ is also a bivector map taking $w_{12}$ to
$w_{21}$.
\end{rem}

\begin{rem}
$R \circ R' = Id_{C_1 \fus C_2} = R' \circ R$
\end{rem}

We now define $R_i : C_1 \fus \cdots\fus (C_i \fus C_{i+1})\fus \cdots
\fus C_n
\to C_1 \fus \cdots \fus (C_{i+1} \fus C_i) \fus \cdots \fus C_n$ to
be the map given by
$$
R_i(g_1, ..., g_i, g_{i+1}, ... g_n) = (g_1, ... , Ad_{g_i}g_{i+1},
g_i, ... , g_n)
$$
that is, $R$ applied to the $i$th and $(i+1)$th term of $M_r$.  $R_i'$
can be defined in a similar way.

\begin{lem}
The full braid group $B_n$ has a faithful representation as a group of
automorphism of the closed $n$-gons in $\s^3$
in which side-lengths are fixed but the order of the sides is not
fixed.  The generators of $B_n$ are given by $R_i$, $1 \leq i \leq
n-1$.
\end{lem}

We now restrict $B_n$ to $P_n$ to get an action of the pure braid
group on $C_r^n$.  This action induces a symplectomorphism on the
moduli space $M_r$.

\begin{cor}
Let $A_{ij} = R_{j-1} \circ \cdots \circ  R_{i+1}
\circ R_i^2 \circ R_{i+1}' \circ \cdots \circ R_{j-1}', \; 1 \leq
i < j \leq n$. $A_{ij}$ induces a symplectomorphism from $M_r$ to
itself. $A_{ij}, \; 1 \leq i < j \leq n$ are the generators of $P_n$
which has a faithful representation as a group of automorphisms of $M_r$.
\end{cor}

We will now show that the braid group actions $A_{ij}$ can be
realized as the time one Hamiltonian flows of the Hamiltonians
$h_{ij}$ given at the start of the section. We begin by studying
the Hamiltonian flows associated to the functions $f_{ij}\in
C^\8(C_r^n)^K$ given by $f_{ij}(g) = tr(g_i g_j)$.  Define
$F_{ij}:C_r^n \to \k$ by $F_{ij}(g) = \big((g_i g_j) - (g_i
g_j)^{-1}\big)$.

The Hamiltonian flow associated to $f_{ij}$ is given by $\Phi_{ij}^t(g)
= (\ha{g_1}(t), ..., \ha{g_n}(t))$ where
\begin{equation*} \ha{g_k}(t) =
\begin{cases}
g_k, \: 0 < k < i \; \textnormal{and} \; j< k < n+1 \\
Ad\Big(\exp\big(t F_{ij}(g)\big)\Big)g_k, \; k = i , j \\
Ad\Big(\exp\big(t F_{ij}(g)\big) g_j \exp\big(- t F_{ij}(g)\big)
g_j^{-1}\Big) g_k, \; i < k < j.
\end{cases}
\end{equation*}

The following formula is used to relate $\Phi_{ij}^t$ to $A_{ij}$.

\begin{lem}
$$
\exp{\left(\frac{\cos^{-1}(-\half tr(g))}{\sqrt{4-tr^2(g)}}
(g-g^{-1})\right)} = g
$$
\end{lem}

We now notice that for time $t = \frac{\cos^{-1}(-\half
f_{ij}(g))}{\sqrt{4 - f_{ij}^2 (g)}}$,

$$
\Phi^t_{ij} = A_{ij}.
$$

The time for which the $\Phi^t_{ij}$ flows
depends on the point in $M_r$ at which flow begins.  We would like
time to be independent on the starting point.  We can achieve this by taking
the Hamiltonian flows of the functions $h_{ij} = \half \left(\cos^{-1}
(-\half f_{ij})\right)^2$.
The Hamiltonian flow $\t{\Phi}_{ij}^t$ associated to
$h_{ij}$ is the renormalization of the flow $\Phi_{ij}^t$ so that
$$
\t{\Phi}_{ij}^1 = A_{ij}
$$
on $M_r$.  We can see the pure braid group as the integer points in
the Hamiltonian flows $\t{\Phi}_{ij}^t, \; 1 \leq i<j \leq n$.

\section{Connection with symplectic forms on relative character
varieties of $n$-punctured 2-spheres}

In this section, we relate the symplectic form on $M_r(\s^3)$ given in
Remark \ref{symform}  to the symplectic form of Goldman type obtained
from the description of $M_r(\s^3)$ as the moduli space of flat
connections on an $n$-punctured 2-sphere. We follow the arguments of
Kapovich and Millson \cite[\S5]{KM1} which considers the analogous question
for $M_r(\E^3)$.  We begin with the general case
in which $G$ is any Lie group with Lie algebra $\g$ which admits a
nondegenerate, $G$-invariant, symmetric, bilinear form.

\subsection{Relative characteristic varieties and parabolic cohomology}

Let $\Si = \s^2-\{p_1, ..., p_n\}$ denote the $n$-punctured 2-sphere
and $U_1, ..., U_n$ be disjoint disc neighborhoods of $p_1, ..., p_n$,
repectively. Further,  $\Ga$ is the fundamental group of $\Si$ with
generators $\ga_i $, $T= \{\Ga_1,
...,\Ga_n\}$ is the collection of subgroups of $\Ga$ with $\Ga_i$  the
cyclic subgroup generated by $\ga_i$, and $U = U_1\cup \cdots \cup U_n$.

Fix $\rho_0 \in \Hom(\Ga, G)$ a representation.  In \cite{KM2} the relative
representation variety $\Hom(\Ga,T;G)$ is defined as the representations
$\rho: \Ga \to G$ such 
that $\rho|_{\Ga_i}$ is contained in the closure of the conjugacy
class of $\rho_0|_{\Ga_i}$.

\begin{rem}
If $G=SU(2)$, there exists a $\rho_0$ such that the relative character 
variety $\Hom(\Ga,T;G)/G$ is isomorphic to $M_r(\s^3)$.  We will make
this isomorphism explicit later on. 
\end{rem}

Let $\rho \in \Hom(\Ga,T;G)$.  Then $\rho$ induces a flat principal
$G$-bundle over $\Si$.  The associated flat Lie algebra bundle will be
denoted by $ad \, P$.  

We define the parabolic cohomology, $H^1_{par} (\Si, ad\; P)$ to
be the subspace of the de Rham cohomology classes in $H^1_{DR} (\Si,
ad\; P)$ whose restrictions to each $U_i$ are trivial.

\subsection{Gauge theoretic description of the symplectic form}

Let $b$ be the nondegenerate, $G$-invariant, symmetric, bilinear
form on $\g$. A skew symmetric bilinear form 
$$
B: H_{par}^1(\Si, ad\, P) \times H_{par}^1(\Si, ad\, P) \to H^2(\Si,
U; \R)
$$ 
is defined  by taking the wedge product together with the bilinear
form $b$.  Evaluating
on the relative fundamental class of $\Si$ gives the skew symmetric form, $$ A:
H_{par}^1(\Si, ad\, P) \times H_{par}^1(\Si, ad\, P) \to \R. $$
Poincare duality give us nondegeneracy of $A$, so $A$ is 
a symplectic form on $\Hom(\Ga,T;G)$.  We will show $A$ corresponds
to the symplectic form $\t\om$ given in Remark \ref{symform}.

We first pass through the group
cohomology description of $H^1_{par} (\Sigma, ad\, P)$ to make this
correspondence explicit.

We identify the universal cover of $\Si$, denoted $\t{\Si}$, with
the hyperbolic plane, $\h^2$.  Let $p:\t\Si \to \Si$ by the covering
projection.  We define the $\mathcal{A}^\bullet(\t\Si , p^* Ad\, P)$ with
$\mathcal{A}^\bullet(\t\Si , \g)$ by parallel translation from a point
$x_0$.  Given $[\eta]\in H^1(\Si,ad\, P)$ choose a representing closed
1-form $\eta\in \mathcal A^1(\Si, ad\, P)$. Let $\t\eta = p^*\eta$.  Then
there is a unique function $f: \t\Si \to \g$ satisfing:
\begin{itemize}
\item $f(x_0) = 0$
\item $df = \t\eta$
\end{itemize}
 
A 1-cochain $h(\eta)\in C^1(\Ga, \g)$ is defined by
$$
h(\eta)(\ga) = f(x) - Ad_\rho(\ga) f(\ga^{-1}x).
$$
This induces an isomorphism from  $H^1(\Si,ad\, P)$ to $H^1(\Ga, \g)$.  It
can be seen that  $[\eta]\in H^1_{par}(\Si,ad\, P)$ if and only if 
$h(\eta)$ restricted to $\Ga_i$ is exact for all $i$.  That is, there
exists an $x_i \in \g$ such that $h(\eta)(\ga_i^k) = x_i -
Ad_\rho(\ga_i^k) x_i$ for each $\ga_i$ a generator of $\Ga$.

We construct the fundamental domain $\mathcal D$ for $\Ga$ operating on
$\h^2$ as in \cite{KM1}.  Choose $x_0$ on $\Si$ and make cuts along
geodesics from $x_0$ to the cusps.  The resulting fundamental
domain $\mathcal D$ is a geodesic $2n$-gon with vertices $v_1, ...,
v_n$ and cusps $v_1^\8, ... , v_n^\8$ ordered so that as we proceed
clockwise around $\D\mathcal D$ we see $v_1, v_1^\8, ..., v_n,
v_n^\8$.  The generator $\ga_i$ fixes 
$v_i^\8$ and satsfies $\ga_i v_{i+1} = v_i$. Let $e_i$ be the
oriented edge joining $v_i$ to $v_i^\8$ and  $\hat{e}_i$ be the
oriented edge joining $v_i^\8$ to $v_{i+1}$. Then $\ga_i
\hat{e}_{i} = - e_i$.

 Let $\rho \in \Hom(\Ga, T;G)$ and $c, c' \in
T_\rho \left(\Hom(\Ga,T;G)/G \right) \simeq H^1_{par}(\Ga, \g)$ be tangent vectors at $\rho$.  The corresponding
elements in $H^1_{par} (\Si, ad P)$ are denoted  $\al$ and $\al'$.  So
$f: \Si \to \g$ which satisfies $df = \t{\al}$ and $f_i(x_0) = 0$.
Let $f(v_i^\8) = x_i$.  Then  
\begin{eqnarray*}
c(\ga_i) & = & f(x) - Ad_{\rho(\ga_i)} f(\ga_i^{-1} x) \\
	& = & f(v_i^\8) - Ad_{\rho(\ga_i)} f(\ga_i^{-1} v_i^\8) \\
	& = & f(v_i^\8) - Ad_{\rho(\ga_i)} f( v_i^\8) \\
	& = & x_i - Ad_{\rho(\ga_i)} x_i.
\end{eqnarray*}
There is an equivalent formulas for $c', \al',$ and $f'$ with $f'(v_i^\8)
= x'_i$.

Let $B_\bullet(\Ga)$ be the bar resolution of $\Ga$.  Thus
$B_k(\Ga)$ is the free $\Z[\Ga]$-module on the symbols $[\ga_1 |
\ga_2 | \cdots | \ga_k]$ with 
$$ \D[\ga_1 | \ga_2 | \cdots |
\ga_k]  =  \ga_1 [ \ga_2 | \cdots | \ga_k] + \sum_{i=1}^{k-1}
(-1)^i [\ga_1 | \cdots | \ga_i \ga_{i+1} | \cdots | \ga_k ] +
(-1)^k [ \ga_1 | \cdots | \ga_{k-1} ]. $$ 
Let $C_k(\Ga) = B_k(\Ga)
\otimes_{\Z[\Ga]} \Z$ with $\Z[\Ga]$ acting on $\Z$ by the
homomorphism $\eps$ defined by 
$$ \eps(\sum_{i=1}^m a_i \ga_i) =
\sum_{i=1}^m a_i. $$
Then $C_k(\ga)$ is the free abelian group on
the symbols $(\ga_1| \cdots | \ga_k) = [\ga_1 | \ga_2 | \cdots |
\ga_k] \otimes 1$ with 
$$ \D(\ga_1 | \ga_2 | \cdots | \ga_k)  =  (
\ga_2 | \cdots | \ga_k) + \sum_{i=1}^{k-1} (-1)^i (\ga_1 | \cdots
| \ga_i \ga_{i+1} | \cdots | \ga_k ) + (-1)^k ( \ga_1 | \cdots |
\ga_{k-1} ). $$ 
A relative fundamental class $F\in C_2(\Ga)$ is
defined by the property 
$$\D F = \sum_{i=1}^n (\ga_i).$$

Let $[\Ga, \D \Ga] = \sum_{i=2}^n (\ga_1 \cdots \ga_{i-1} | \ga_i)\in
C_2(\Ga)$, then

\begin{lem}
$[\Ga, \D\Ga]$ is a relative fundamental class.
\end{lem}
\proof  The proof is left to the reader.

We will now give the symplectic form $A$ in terms of group
cohomology.  We denote by $\cup_b$ the cup product of
Eilenberg-MacLane cochains using the form $b$ on the coefficients.

\begin{prop}\label{gold}
$$
A(\al, \al') = \sum_{i=1}^n \< c \cup_b x'_i) , (\ga_i) \> - \<
c\cup_b c', [\Ga, \D \Ga] \>
$$
\end{prop}
\no We will use the next Lemmas to prove Proposition \ref{gold}.

\begin{lem}
$$
\int_{e_i} B(f,\t\al') + \int_{\ha{e}_i} B(f,\t\al') = b\left( c(\ga_i),
f'(v_i^\8)\right) - b\left( c(\ga_i), f'(v_i)\right)
$$
\end{lem}
\proof
Recall $\ga_i \ha{e}_i = -e_i$, so that $\ha{e}_i = - \ga_i^{-1} e_i$.
We then have
\begin{eqnarray*}
\int_{e_i} B(f,\t\al') + \int_{\ha{e}_i} B(f,\t\al') & = & \int_{e_i}
B(f,\t\al') + \int_{\ha{e}_i} B(f,\t\al') \\
    & = & \int_{e_i} B(f,\t\al') + \int_{\ga_i^{-1} e_i} B(f,\t\al')\\
    & = & \int_{e_i} B\left(f,\t\al'\right) + \int_{e_i}(\ga_i^{-1})^*
B(f,\t\al') \\
    & = & \int_{e_i} B\left(f,\t\al'\right) + \int_{e_i}
    B\left((\ga_i^{-1})^*f, (\ga_i^{-1})^*\t\al'\right) \\
    & = & \int_{e_i} B\left(f,\t\al'\right) + \int_{e_i}
    B\left(Ad_{\rho(\ga_i)} 
(\ga_i^{-1})^*f, Ad_{\rho(\ga_i)} (\ga_i^{-1})^*\al'\right) \\
    & = &  \int_{e_i} B\left(f - Ad_{\rho(\ga_i)} (\ga_i^{-1})^*f,
\t\al'\right) \\
    & = & \int_{e_i} B\left(c(\ga_i),\t\al'\right) \\
    & = & b\left( c(\ga_i), f'(v_i^\8)\right) - b\left( c(\ga_i),
f'(v_i)\right)
\end{eqnarray*}
\qed

\begin{lem}
$$\sum_{i=1}^n b\left(c(\ga_i), f'(v_i)\right) = \sum_{i=1}^n
b\left(c(\ga_i), f'(v_i^\8)\right) - \sum_{i=1}^n \< c \cup_b y_i,
(\ga_i)\> + \< c\cup_b c', [\Ga, \D \Ga] \>
$$
\end{lem}
\proof
By definition, for any $x \in \h^2$ and $\ga \in \Ga$ we have
$$
c'(\ga) = f'(x) - Ad_{\rho(\ga)} f'(\ga^{-1}x)
$$
Let $\ga = \ga_i$ and $x = v_i$, then
$$
c'(\ga_i) = f'(v_i) - Ad_{\rho(\ga_i)} f'(v_{i+1})
$$
Using  $f'(v_1) = 0$, we obtain
\begin{eqnarray*}
c'(\ga_1 \cdots \ga_i) & = & f'(v_1) - Ad_{\rho(\ga_1 \cdots \ga_i)}
    f'(\ga_i^{-1} \cdots \ga_1^{-1} v_{1}) \\
    & = &  - Ad_{\rho(\ga_1 \cdots \ga_i)} f'(v_{i+1}).
\end{eqnarray*}

\no We will also need
\begin{eqnarray*}
c'(\ga_1 \cdots \ga_i)  & = & c'(\ga_1 \cdots \ga_{i-1}) + Ad_{\rho(\ga_1
\cdots \ga_{i-1})} c'(\ga_i) \\
    & = & c'(\ga_1) + Ad_{\rho(\ga_1)}c'(\ga_2) + \cdots +  Ad_{\rho(\ga_1
\cdots \ga_{i-1})} c'(\ga_i)
\end{eqnarray*}
\no and, since $\ga_1 \cdots \ga_n = 1$,
$$
0 = c'(\ga_1 \cdots \ga_n) = c'(\ga_1) + Ad_{\rho(\ga_1)}c'(\ga_2) +
\cdots +  Ad_{\rho(\ga_1 \cdots \ga_{n-1})} c'(\ga_n)
$$

We then have,

\begin{eqnarray*}
\sum_{i=1}^n b\left(c(\ga_i), f'(v_i)\right) & = & - \sum_{i=1}^n
b\left(c(\ga_i),  Ad_{\rho(\ga_1 \cdots \ga_i)^{-1}}c'(\ga_1 \cdots
\ga_{i-1})\right) \\
     & = & - \sum_{i=1}^n b\left(Ad_{\rho(\ga_1
\cdots \ga_{i-1})} c(\ga_i), c'(\ga_1) +
Ad_{\rho(\ga_1)}c'(\ga_2) + \cdots +  Ad_{\rho(\ga_1 \cdots
\ga_{i-2})} c'(\ga_{i-1})\right) \\
     & = & - \sum_{i=1}^n \sum_{j=1}^{i-1} b\left(Ad_{\rho(\ga_1
\cdots \ga_{i-1})} c(\ga_i),  Ad_{\rho(\ga_1 \cdots
\ga_{j-1})} c'(\ga_{j})\right) \\
     & = & - \sum_{j=1}^n \sum_{i=j+1}^{n} b\left(Ad_{\rho(\ga_1
\cdots \ga_{i-1})} c(\ga_i),  Ad_{\rho(\ga_1 \cdots
\ga_{j-1})} c'(\ga_{j})\right) \\
    & = & \sum_{j=1}^n \sum_{i=1}^{j} b\left(Ad_{\rho(\ga_1
\cdots \ga_{i-1})} c(\ga_i),  Ad_{\rho(\ga_1 \cdots
\ga_{j-1})} c'(\ga_{j})\right) \\
    & = & \sum_{j=1}^n  b\left( c(\ga_1 \cdots \ga_j),
Ad_{\rho(\ga_1 \cdots \ga_{j-1})} c'(\ga_{j})\right) \\
    & = & \sum_{j=1}^n  b\left( c(\ga_1 \cdots \ga_{j-1}) + Ad_{\rho(\ga_1
\cdots \ga_{j-1})} c(\ga_j) , Ad_{\rho(\ga_1 \cdots \ga_{j-1})}
c'(\ga_{j})\right) \\
    & = & \sum_{j=1}^n  b\left( c(\ga_1 \cdots \ga_{j-1}) ,
Ad_{\rho(\ga_1 \cdots \ga_{j-1})} c'(\ga_{j})\right) + \sum_{j=1}^n
b\left( c(\ga_j) , c'(\ga_{j})\right) \\
    & = & \<c \cup_b c'), [\Ga,\D \Ga]\> +  \sum_{j=1}^n b\left( c(\ga_j) ,
f'(v^\8_{j}) - Ad_{\rho(\ga_j)}f'(v_j^\8)\right) \\
    & = & \<c\cup_b c', [\Ga,\D \Ga]\> +  \sum_{j=1}^n b\left( c(\ga_j) ,
f'(v^\8_{j})\right) -  \sum_{j=1}^n \<B(c,y'_j), (\ga_j)\>
\end{eqnarray*}
\qed

\no{\it Proof of Proposition \ref{gold}:}
\begin{eqnarray*}
A(\al, \al') & = & \int_{\Sigma} B(\al, \al') \\
	& = & \int_{\mathcal D} B(\t\al, \t\al') \\
	& = & \int_{\D\mathcal D} B(\t\al, f') \\
	& = & \sum_{i=1}^n\left( \int_{e_i} B(\t\al, f') +
\int_{\ha{e}_i} B(\t\al, f')\right) \\ 
    & = & \sum_{j=1}^n \<c\cup_b x'_j), (\ga_j)\> -  \< c\cup_b c',
[\Ga,\D \Ga]\>
\end{eqnarray*}\qed

\subsection{Correspondence between $M_r(\s^3)$ and
$\Hom\left(\Ga,T;SU(2)\right)/SU(2)$}

We now restrict to the case $G = SU(2)$.  We define the isomorphism $$
\Ups:\Hom\left(\Ga,T;SU(2)\right) \to \t{M}_r, $$ 
where $\t{M}_r$ is the closed polygonal linkages in $\s^3$ based
at a point, by $$ \Ups(\rho) = \left( \rho(\ga_1), ... ,
\rho(\ga_n)\right). $$ This induces an isomorphism, which we also
denote by $\Ups$, $$ \Ups:\Hom(\Ga,T;SU(2))/SU(2) \to M_r. $$

\no The differential $d\Ups_\rho: T_\rho (\Hom(\Ga,T;SU(2))/SU(2)) \to
T_{\Ups(\rho)}M_r$ is then defined by $$d\Ups_\rho(c) = \left(
dR_{\rho(\ga_1)}c(\ga_1), ... , dR_{\rho(\ga_n)}c(\ga_n)\right).
$$ Here $T_\rho (\Hom(\Ga,T;SU(2))/SU(2))$ is identified with an element of $\Z^1_{par}(\Ga,\g)$. We have $$
 d\Ups_\rho(c) = \left( dR_{g_1} x_1 - dL_{g_1} x_1, ... ,
dR_{g_n} x_n - dL_{g_n} x_n \right) $$ and $$
 d\Ups_\rho(c') = \left( dR_{g_1} x'_1 - dL_{g_1} x'_1, ... ,
dR_{g_n} x'_n - dL_{g_n} x'_n \right). $$

Recall, the symplectic form on $M_r$ is given by $$ \t\om = \sum_{i=1}^n  \om_i + \half \sum_{i=1}^n
\sum_{j=i+1}^n \left( Ad_{g_1 \cdots g_{i-1}} \bar\th_i \wedge_b
Ad_{g_1 \cdots g_{j-1}} \bar\th_j \right). $$ We can now prove the
main result of this section
\begin{thm}
$\Ups^*\t\om = A$
\end{thm}
\proof

First we note that
$$
\Ups^*\bar\th_i(c) = c(\ga_i)
$$
and
\begin{eqnarray*}
(\Ups^* \om_i)(c,c') & = & \om_i \left(dR_{g_i} c(\ga_i), dR_{g_i}
c'(\ga_i)\right) \\
    & = & - \half \left( Ad_{g_i^{-1}} c(\ga_i) + c(\ga_i), x'_i
\right) \\
    & = & - \half \left( c(\ga_i),  Ad_{g_i} x'_i + x'_i
\right) \\
    & = & - \half \left( c(\ga_i), c'(\ga_i) \right) - \left
( c(\ga_i),  Ad_{g_i} x'_i \right) \\
    & = & - \half \left( Ad_{g_1 \cdots g_{i-1}} c(\ga_i),
Ad_{g_1 \cdots g_{i-1}} c'(\ga_i) \right) +  \<c \cup_b x'_i), (\ga_i)\> \\
\end{eqnarray*}
It follows that

\begin{eqnarray*}
(\Ups^*\t\om)(c,c') & = & \sum_{i=1}^n (\Ups^* \om_i)(c,c') + \half
\sum_{i=1}^n   \sum_{j=i+1}^n \Ups^* \left( Ad_{g_1 \cdots g_{i-1}}
\bar\th_i \wedge_b Ad_{g_1 \cdots g_{j-1}} \bar\th_j \right)(c,c') \\
    & = & \sum_{i=1}^n \<c \cup_b x'_i), (\ga_i)\>  - \sum_{i=1}^n
\half \left( Ad_{g_1 \cdots g_{i-1}} c(\ga_i), Ad_{g_1 \cdots g_{i-1}}
c'(\ga_i) \right) \\
    &   & + \sum_{i=1}^n \sum_{j=i+1}^n \left( Ad_{g_1
\cdots g_{i-1}} c(\ga_i) , Ad_{g_1 \cdots g_{j-1}} c'(\ga_j) \right)
\\
    &   & - \sum_{i=1}^n \sum_{j=i+1}^n \left( Ad_{g_1
\cdots g_{i-1}} c'(\ga_i) , Ad_{g_1 \cdots g_{j-1}} c(\ga_j) \right)\\
    & = & \sum_{i=1}^n \< c \cup_b x'_i , (\ga_i)\>  - \sum_{i=1}^n
\half \left( Ad_{g_1 \cdots g_{i-1}} c(\ga_i), Ad_{g_1 \cdots g_{i-1}}
c'(\ga_i) \right) \\
    &   & + \sum_{j=2}^n \sum_{i=1}^{j-1} \left( Ad_{g_1
\cdots g_{i-1}} c(\ga_i) , Ad_{g_1 \cdots g_{j-1}} c'(\ga_j) \right)
\\
    &   & + \sum_{i=1}^n \sum_{j=1}^i \left( Ad_{g_1
\cdots g_{i-1}} c'(\ga_i) , Ad_{g_1 \cdots g_{j-1}} c(\ga_j) \right)\\
    & = & \sum_{i=1}^n \<c\cup_b x'_i, (\ga_i)\> + \sum_{j=2}^n
\sum_{i=1}^{j-1} \left( Ad_{g_1 \cdots g_{i-1}} c(\ga_i) ,
Ad_{g_1 \cdots g_{j-1}} c'(\ga_j) \right) \\
    & = & \sum_{i=1}^n \<c \cup_b x'_i), (\ga_i)\> + \sum_{j=2}^n
\left( Ad_{g_1 \cdots g_{i-1}} c'(\ga_i) ,  c(\ga_1 \cdots \ga_{i-1})
\right) \\
    & = & \sum_{i=1}^n \< c \cup_b x'_i, (\ga_i)\> -  \< c\cup_b c',
[\Ga,\D \Ga]\> \\
	& = & A(\al, \al')
\end{eqnarray*}
\qed

It is easily seen that the functions $\ell_i$ from \S \ref{commflow}
corresponds to the following Goldman functions.  Let $\varphi : G \to
\R$ be defined by $\varphi(g) = \cos^{-1}\left(-\half trace(g) \right)$.
We then defined the function $\varphi_\ga :
\Hom\left(\Ga,T;SU(2)\right)/SU(2) \to \R$ by $\varphi_ga(\rho) =
\varphi\left(\rho(-ga)\right)$.  We see that $$\Ups^*\ell_i =
\varphi_{\ga_1 \cdots \ga_i}$$ Then choosing an maximal collection of
nonintersecting diagonal on $M_r$ corresponds to a pair of pants
decomposition on $\Si$.


\begin{thebibliography}{99999}
\addcontentsline{toc}{section}{Bibliography}


\bibitem[AKS]{AKS}
A.\ Alekseev, Y. Kosmann-Schwarzbach, {\em Manin pairs and moment maps },
preprint, math.DG/9909176.

\bibitem[AKSM]{AKSM}
A.\ Alekseev, Y. Kosmann-Schwarzbach, E. Meinrenken, {\em
Quasi-Poisson Manifolds}, prepint math.DG/0006168.

\bibitem[AMM1]{AMM}
A.\ Alekseev, A. Malkin, E. Meinrenken, {\em Lie group valued
moment maps}, J. Differential Geom. {\bf 48} (1998), 445--495.

\bibitem[AMM2]{AMM2}
A.\ Alekseev, A. Malkin, E. Meinrenken, {\em Manin pairs of a compact
simple Lie algebra}, unpublished notes.

\bibitem[Bi]{Bi}
J.\ Birman, ``Braids, links, and mapping class groups'', Annals of
Mathematics Studies, No. 82,  Princeton Univ. Press, 1974.


\bibitem[CP]{CP}
V.\ Chari, A.\ Pressley, ``A guide to quantum groups'', Cambridge
Univ. Press, 1994.


\bibitem[FM]{FM}
H.\ Flaschka, J.J.\ Millson, {\em An integrable system on the
moduli space of n points in $\mathbb{CP}^m$}, in preparation.

\bibitem[FR]{FR}
H.\ Flaschka, T.\ Ratiu,  {\em A convexity theorem for Poisson
actions of compact Lie groups}, Ann. Sci. Ecole Norm. Sup., (4)
vol. {\bf 29} (1996), no. 6, 787--809.

\bibitem[Go]{Go}
W.\ Goldman, {\em Invariant functions on Lie groups and Hamiltonian
flows of surface group representations}, Invent. math. {\bf 85}
(1986), 263--302.



\bibitem[GHJW]{GHJW} K.\ Guruprasad, J.\ Huebshmann, L.\ Jeffrey, A.\
Weinstein, {\em Group systems, groupoids, and moduli spaces of
parabolic bundles},  Duke Math. J. {\bf 89} (1997), 377--412.

\bibitem[Je]{Je}
L.\ Jeffrey, {\em Extended moduli spaces of flat connections on
Riemann surfaces}, Math. Ann. {\bf 298} (1994), 667--692.


\bibitem[KLM]{KLM}
M.\ Kapovich, B.\ Leeb, J.J.\ Millson, in preparation.

\bibitem[KM1]{KM1}
M.\ Kapovich, J.J.\ Millson,  {\em The symplectic geometry of
polygons in Euclidean space},  J. Differential Geom. {\bf 44}
(1996), 479--513.

\bibitem[KM2]{KM2}
M.\ Kapovich, J.J.\ Millson,  {\em The relative deformation theory of
representations of flat connections and deformations of linkages in
constant curvature spaces},  Compositio Mathematica {\bf 103} (1996),
287--317.

\bibitem[KMT]{KMT}
M.\ Kapovich, J.J.\ Millson, T.\ Treloar  {\em The symplectic
geometry of polygons in hyperbolic 3-space}, Asian Journal
Math {\bf 4} (2000), 123--164.


\bibitem[Ki]{Ki}
F.\ Kirwan,  {\em Cohomology of quotients in symplectic and
algebraic geometry},  Mathematical Notes, Princeton University
Press, 1984.

\bibitem[KS1]{KS1}
Y.\ Kosmann-Schwarzbach,  {\em Lie bialgebras, Poisson Lie groups, and
dressing transformations},  Lecture Notes Phys.  {\bf 495}, Springer-Verlag,
1997, 104--170.

\bibitem[KS2]{KS2}
Y.\ Kosmann-Schwarzbach,  {\em Jacobian quasi-bialgebras and
quasi-Poisson Lie groups},  Contemporary Mathematics {\bf 132} (1991),
459--489.

\bibitem[LM]{LM}
B.\ Leeb, J.J.\ Millson,  {\em Convex functions on symmetric
spaces and geometric invariant theory for weighted configurations
in flag manifolds}, in preparation.

\bibitem[Le]{Le}
M.\ Leingang, {\em Symmetric pairs and moment spaces}, preprint
math.SG/9810064.

\bibitem[Lu1]{Lu}
J.-H.\ Lu, {\em Multiplicative and affine Poisson structures on
Lie groups}, Ph.D. thesis, University of California, Berkeley,
1990.




\bibitem[Mi]{Mi}
J.J.\ Millson, {\em Bending polygons and decomposing tensor
products, three examples}, in preparation.


\bibitem[MZ]{MZ}
J.J.\ Millson, B.\ Zombro, {\em A K\"{a}hler structure on the
moduli space of isometric maps of a circle into Euclidean space},
Invent. Math., {\bf 123} (1996), 35--59.


\bibitem[STS]{Se}
M.\ Semenov-Tian-Shansky, {\em Dressing transformations and
Poisson  group actions}, Publ. of RIMS, vol. {\bf 21} (1985)
1237--1260.

\bibitem[Sh]{Sh}
R.\ W.\ Sharpe, ``Differential geometry'', Graduate Texts in Math.,
vol. 166, Springer, 1997.

\bibitem[Th]{Th}
W.\ Thurston, ``Three-dimensional geometry and topology'', Princeton
Univ. Press, 1997 .

\bibitem[Tr]{Tr}
T.\ Treloar, {\em The symplectic geometry on loops in the
3-sphere}, in preparation.
\end{thebibliography}
\end{document}